% Logic Eprints
%Submitted 0950 Wed Sep 22, 1993 by: andreas.blass@math.lsa.umich.edu (andreas blass )
%logic/blass/quest-ans.tex
%

%This file is for use with AMSTeX (version 2.1).
\magnification=1200
\input amstex
\documentstyle{amsppt}

%macros go here
%The next few lines are Vaughan Pratt's macros for drawing a par
%symbol (an inverted ampersand) called \invamp or (boxed) \Invamp
\newcount\PLv\newcount\PLw\newcount\PLx\newcount\PLy\newdimen\PLyy\newdimen\PLX
\newbox\PLdot \setbox\PLdot\hbox{.} \def\scl{.08} % resettable scale
\def\PLot#1{\PLx`#1\advance\PLx-42\PLy\PLx\PLv\PLx\divide\PLy9\PLw\PLy\multiply
\PLw9\advance\PLx-\PLw\advance\PLx-4\PLy-\PLy\advance\PLy4\PLX=\the\PLx
pt
\advance\PLyy\the\PLy pt\wd\PLdot=\scl\PLX\raise\scl\PLyy\copy\PLdot}
\def\draw#1{\ifx#1\end\let\next=\relax\else\PLot#1\let\next=\draw\fi\next}

%  Girard's inverted ampersand.  Usage: \invamp. Draw time @ 10 Mips:
%  .33 sec
\def\invamp{\hbox{\PLyy=70pt\draw
:::;DMV_gqppyyyyyooooxxxnnwvlutkjaWNE=5-./9%
9:::CCCC:::99/..--544=EENWWaajjjkktttttttNNNVVVVVVVV\end \hskip4pt}}
%  \Invamp = Boxed \invamp.  Draw time < .02 sec, but max ~300
%  chars/page
\newbox\iabox\setbox\iabox\invamp \def\Invamp{\copy\iabox}
\define\pa{\Invamp\,}
\define\dpa{\overline{\pa}}
\define\ti{\otimes}
\define\dti{\overline\otimes}
\define\pv{$\Cal{PV}$}
\define\ba{\bold{A}}
\define\bbb{\bold{B}}
\define\obj#1{(#1_-,#1_+,#1)}
\define\cc{\frak{c}}
\define\rr{\frak{r}}
\define\rs{{\frak r}_\sigma}
\redefine\ss{\frak{s}}
\define\dd{\frak{d}}
\define\bb{\frak{b}}
\define\nn{\Bbb{N}}
\define\fn{{\nn^\nn}}
\define\pn{\Cal P(\nn)}
\define\pin{\Cal P_\infty(\nn)}
\define\imp{\multimap}
\define\voj{Vojt\'a\v s}

\topmatter
\title
Questions and Answers --- A Category Arising in \\
Linear Logic, Complexity Theory, and Set Theory
\endtitle
\rightheadtext{Questions and Answers}
\author
Andreas Blass
\endauthor
\address
Mathematics Dept., University of Michigan, Ann Arbor, 
MI 48109, U.S.A.
\endaddress
\email
ablass\@umich.edu
\endemail
\thanks Partially supported by NSF grant DMS-9204276.
\endthanks
\abstract
A category used by de~Paiva to model linear logic also occurs in
Vojt\'a\v s's analysis of cardinal characteristics of the continuum.
Its morphisms have been used in describing reductions between search
problems in complexity theory.  We describe this category and how it
arises in these various contexts.  We also show how these contexts
suggest certain new multiplicative connectives for linear logic.
Perhaps the most interesting of these is a sequential composition
suggested by the set-theoretic application.
\endabstract
\endtopmatter
\document

\head
Introduction
\endhead

The purpose of this paper is to discuss a category that has appeared
explicitly in work of de~Paiva \cite{15} on linear logic and in work of
\voj\ \cite{21, 22} on cardinal characteristics of the continuum.  We call
this category \pv\ in honor of de~Paiva and \voj\ (or, more informally,
in honor of Peter and Valeria).  The same category is implicit in a
concept of many-one reduction of search problems in complexity theory
\cite{12, 19}. 

The objects of \pv\ are binary relations between sets; more precisely
they are triples $\ba=\obj A$, where $A_-$ and $A_+$ are sets and $A
\subseteq A_-\times A_+$ is a binary relation between them.  (We
systematically use the notation of boldface capital letters for
objects, the corresponding lightface letters for the relation
components, and subscripts $-$ and $+$ for the two set components.) A
morphism from $\ba$ to $\bbb=\obj B$ is a pair of functions
$f_-:B_-\to A_-$ and $f_+:A_+\to B_+$ such that, for all $b\in B_-$
and all $a\in A_+$,
$$
A(f_-(b),a)\implies B(b,f_+(a)).
$$
(Note that the function with the minus subscript goes backward.)
Composition of these morphisms is defined componentwise, with the
order reversed on the minus components: $(f\circ g)_-=g_-\circ f_-$
and $(f\circ g)_+=f_+\circ g_+$.  This clearly defines a category \pv.

The category \pv\ is the special case of de~Paiva's construction
$\bold{GC}$ from \cite{15} where $\bold C$ is the category of sets.  It
is also the dual of \voj's category $GT$ of generalized Galois-Tukey
connections \cite{21, 22}.

Intuitively, we think of an object $\ba$ of \pv\ as representing a
problem (or a type of problem).  The elements of $A_-$ are instances
of the problem, i.e., specific questions of this type; the elements
of $A_+$ are possible answers; and the relation $A$ represents
correctness, i.e., $A(x,y)$ means that $y$ is a correct answer to the
question $x$.

There are strong but superficial similarities between \pv\ and a
special case of a construction due to Chu and presented in the
appendix of \cite1 and Section~3 of \cite2.  (Readers unfamiliar with
the Chu construction can skip this paragraph, as it will not be
mentioned later.)  Specifically, Chu's construction, applied to the
cartesian closed category of sets and the object 2, yields a
$*$-autonomous category in which the objects are the same as those of
\pv\ and the morphisms differ from those of \pv\ only in that they are
required to satisfy $A(f_-(b),a)\iff B(b,f_+(a))$ rather than just an
implication from left to right.  This apparently minor difference in
the definition leads to major differences in other aspects of the
category.  Specifically, the internal hom-functor and the tensor
product in Chu's category are entirely different from those of \pv.

In the next few sections, we shall describe how \pv\ arose in various
contexts. Thereafter, we indicate how ideas that arise naturally in
these contexts suggest new constructions in linear logic.

\head
Reductions of Search Problems
\endhead

Much of the theory of computational complexity (e.g., \cite8) deals
with decision problems.  Such a problem is specified by giving a
set of instances together with a subset called the set of positive
instances; the problem is to determine, given an arbitrary instance,
whether it is positive.  In a typical example, the instances might be
graphs and the positive instances might be the 3-colorable graphs.  In
another example, instances might be boolean formulas and positive
instances might be the satisfiable ones.  A {\sl (many-one)
reduction\/} from one decision problem to another is a map sending
instances of the former to instances of the latter in such a way that
an instance of the former is positive if and only if its image is
positive.  Clearly, an algorithm computing such a reduction and an
algorithm solving the latter decision problem can be combined to yield
an algorithm solving the former.

There are situations in complexity theory where it is useful to
consider not only decision problems but also search problems.  A
search problem is specified by giving a set of instances, a set of
witnesses, and a binary relation between them; the problem is to find,
given an instance, some witness related to it.  For example, the
3-colorability decision problem mentioned above (given a graph, is it
3-colorable?) can be converted into the 3-coloring search problem
(given a graph, find a 3-coloring).  Here the instances are graphs,
the witnesses are 3-valued functions on the vertices of graphs, and
the binary relation relates each graph to its (proper) 3-colorings.
Similarly, there is a search version of the boolean satisfiability
problem, where instances are boolean formulas, witnesses are truth
assignments, and the binary relation is the satisfaction relation.
Notice that a search problem is just an object $\ba$ of \pv, the set of
instances being $A_-$ and the set of witnesses $A_+$.  

There is a reasonable analog of many-one reducibility in the context
of search problems.  A reduction of $\bbb$ to $\ba$ should first
convert every instance $b\in B_-$ of $\bbb$ to an instance $a\in A_-$
of $\ba$ (just as for decision problems), and then, if a witness $w$
related to $a$ is given, it should allow us, using $w$ and remembering
the original instance $b$, to compute a witness related to $b$.
Again, an algorithm computing such a reduction and an algorithm
solving $\ba$ can clearly be combined to yield an algorithm solving
$\bbb$.  Most known many-one reductions between NP decision problems
\cite8 implicitly involve many-one reductions of the corresponding
search problems.

Formally, a {\sl reduction\/} therefore consists of two functions,
$f_-:B_-\to A_-$ and $f_+:A_+\times B_-\to B_+$ such that, for all
$b\in B_-$ and $w\in A_+$,
$$
A(f_-(b),w)\implies B(b,f_+(w,b)).
$$
This is nearly, but not quite, the definition of a morphism from $\ba$
to $\bbb$.  The difference is that in a morphism $f_+$ would have only
$w$, not $b$, as its argument.  Thus, morphisms amount to reductions
where the final witness (for $b$) is computed from a witness $w$ for
$a=f_-(b)$ without remembering $b$.  This notion of reduction has been
used in the literature \cite{12, 19}, but I would not argue that it is as
natural as the version where one is allowed to remember $b$.

These observations lead to a suggestion that we record for future
reference. 
\proclaim{Suggestion 1}
Find a natural place in the theory of \pv\ for reductions as described
above, i.e., pairs of functions that are like morphisms except that
$f_+$ takes an additional argument from $B_-$ and the implication
relating $f_-$ and $f_+$ is amended accordingly.
\endproclaim

A ``dual'' modification of the notion of morphism, allowing $f_-$ to
have an extra argument in $A_+$, occurred in de~Paiva's work \cite{14} on
a categorial version of G\"odel's Dialectica interpretation, work that
preceded the introduction of \pv\ in \cite{15}.

\head
Linear Logic
\endhead

The search problems (objects of \pv) and reductions (morphisms of \pv\
or generalized morphisms as in Suggestion~1) described in the
preceding section are vaguely related to some of the intuitions that
underlie Girard's linear logic \cite9.  Girard has written about
linear logic as a logic of questions and answers (or actions and
reactions) \cite{9, 10}, so it seems reasonable to try to model this
idea in terms of \pv.  Also, the fact that in a many-one reduction of
$\bbb$ to $\ba$ a witness for $\bbb$ is produced from exactly one
witness for $\ba$ is reminiscent of the central idea of linear logic
that a conclusion is obtained by using each hypothesis exactly once.
In this section, we attempt to make these vague intuitions precise.
Our goal here is to develop de~Paiva's interpretation of linear logic
(at least the multiplicative and additive parts; the exponentials will
be discussed briefly later) in a step by step fashion that emphasizes
the naturality or necessity of the definitions used.

We intend to use objects of \pv\ as the interpretations of the
formulas of linear logic.  This corresponds to Girard's intuition that
for any formula $A$ there are questions and answers of type $A$.  Of
course, in addition to questions and answers, objects of \pv\ also
have a correctness relation between them.  It is reasonable to expect
that one formula linearly implies another, in  a particular
interpretation, if and only if there is a morphism in \pv\ from (the
object interpreting) the former to (the object interpreting) the
latter; we shall see this more precisely later.  

To produce an interpretation of linear logic, we must tell how to
interpret the connectives, and we must define what it means for a
sequent to be true in an interpretation.

Perhaps the easiest part of this task is to interpret the additive
connectives, $\&$ and $\oplus$.  It seems to be universally accepted
\cite{17} 
that a reasonable categorial model of linear logic will interpret
these as the product and coproduct of the category.  Fortunately, \pv\
has products and coproducts, so we adopt these as the interpretations
of the additive connectives.  The result is that ``with'' is
interpreted as
$$
\obj A\&\obj B=(A_-+B_-,A_+\times B_+,W),
$$
where
$$
W(x,(a,b))\iff\cases
A(x,a),&\text{if }x\in A_-\\
B(x,b),&\text{if }x\in B_-;
\endcases
$$
``plus'' is interpreted as
$$
\obj A\oplus\obj B=(A_-\times B_-,A_++B_+,V)
$$
where
$$
V((a,b),x)\iff\cases
A(a,x),&\text{if }x\in A_+\\
B(b,x),&\text{if }x\in B_+;
\endcases
$$
and the additive units are 
$$
\top=(\emptyset, 1, \emptyset)\quad\text{and}\quad
0=(1,\emptyset, \emptyset),
$$
where 1 represents any one-element set.

These definitions correspond reasonably well to the intuitive meanings
of the additive connectives in terms of questions and answers or in
terms of Girard's ``action'' description of linear logic \cite{10}.  To
answer a disjunction $A\oplus B$ is to provide an answer to one of $A$
and $B$; correctness means that, confronted with questions of both
types, we answer one of them correctly (in the sense of $A$ or $B$).
To answer a conjunction $A\& B$ we must give answers for both, but we
are confronted with a question of only one type and only our answer to
that one needs to be correct.  The intuitive discussion of
conjunction, in particular the fact that we must give answers of both
types even though only one will be relevant to the question, might
make better sense if we think of the answer as being given before the
question is known.  This is a rather strange way of running a
dialogue, but it will arise again later in other contexts (and I've
seen examples of it in real life).

There is also a natural interpretation of linear negation, since (cf.
\cite{9, 10})  questions of type $A$ are answers of type the negation
$A^\perp$ of $A$ and vice versa.  We define 
$$
{\obj A}^\perp=(A_+,A_-,A^\perp), 
$$
where
$$
A^\perp(x,y)\iff\neg A(y,x).
$$
So linear negation interchanges questions with answers and replaces
the correctness relation by the complement of its converse.  Perhaps a
few words should be said about the use of the complement of the
converse rather than just the converse.  There are several reasons for
this, perhaps the most intuitive being that we are, after all,
defining a sort of negation.  Another way to look at it is to think of
a contest between a questioner and an answerer, where success for the
questioner is defined to mean failure for the answerer (cf. the
discussion of challengers and solvers in \cite{11}).  ``That's a good
question'' often means that I have no good answer.  For another
indication that the given definition of ${}^\perp$ is appropriate, see
the section on set-theoretic applications below.  

Mathematically, the strongest reason for defining ${}^\perp$ as we did
is that it gives a contravariant involution of the category \pv.  That
is, the operation ${}^\perp$ on objects and the operation on morphisms
defined by $(f_-,f_+)^\perp=(f_+,f_-)$ constitute a contravariant
functor from \pv\ to itself, whose square is the identity.  This
corresponds to the equivalences in linear logic between $A\vdash B$ and
$B^\perp\vdash A^\perp$ and between $A^{\perp\perp}$ and $A$.

We turn now to a more delicate matter, the interpretation of the
multiplicative connectives.  We begin with ``times.''  Girard's
intuitive explanation of the difference between the multiplicative
conjunction $\ti$ and the additive conjunction $\&$ in \cite{10} is that
the former represents an ability to perform both actions while the
latter represents an ability to do either one of the two actions
(chosen externally).  Looking back at the interpretation of $\&$, we
would expect to modify it by allowing questions of both sorts, rather
than just one, and requiring both components of the answer to be
correct.  This operation on objects of \pv\ is quite natural, and
occurs in both \cite{15} and \cite{21}.  De~Paiva uses the notation $\ti$
for it, although it is not the interpretation of Girard's connective
$\ti$ in her interpretation of linear logic.  \voj\ uses the notation
$\times$ even though it is not the product in the category.  We shall
use the notation $\dti$ and regard it as a sort of provisional tensor
product.  Formally, we define
$$
\obj A\dti\obj B=(A_-\times B_-,A_+\times B_+, A\times B),
$$
where the relation $A\times B$ is defined by
$$
(A\times B)((x,y),(a,b))\iff A(x,a)\text{ and }B(y,b).
$$
Of course, since we have already interpreted negation, our provisional
$\ti$ gives rise to a dual connective, the provisional ``par'':
$$
\obj A\dpa\obj B=(A_-\times B_-,A_+\times B_+,P)
$$
where
$$
P((x,y),(a,b))\iff A(x,a)\text{ or }B(y,b).
$$
To see why these interpretations of the multiplicative connectives are
only provisional and must be modified, we turn to the question of
soundness of the interpretation.  This requires, of course, that we
define what is meant by a sequent being valid, which presumably
depends on a notion of sequents being true in particular
interpretations, i.e, with particular objects as values of the atomic
formulas.  For simplicity, we work with one-sided sequents, as
in \cite9.  So a sequent is a finite list (or multi-set) of formulas,
each interpreted as an object of \pv.  Since a sequent is deductively
equivalent in linear logic with the par of its members, we interpret
the sequent as the (provisional) par of its members, i.e., as a
certain object of \pv.  So we must specify what we mean by truth of
an object of \pv, and then we must try to verify the soundness of the
axioms and rules of linear logic.

There are two plausible interpretations of truth of an object
$\ba=\obj A$, both saying intuitively that one can answer all the
questions of type $\ba$.  The difference between the two is in whether
the answer can depend on the question.  

The first (provisional) interpretation of truth allows the answer
to depend on the question, as one would probably expect intuitively.
$$
\models_1\obj A\iff\forall x\in A_-\,\exists y\in A_+\,A(x,y).
$$
The second, stronger (provisional) interpretation is that one answer
must uniformly answer all questions correctly.
$$
\models_2\obj A\iff\exists y\in A_+\,\forall x\in A_-\,A(x,y).
$$
Before dismissing the second interpretation as unreasonably strong,
one should note that the two interpretations are dual to each other in
the sense that $\ba$ is true in either sense if and only if its
negation $\ba^\perp$ is not true in the other sense.  Furthermore, the
second definition fits better with the idea that truth of a sequent
$A\vdash B$ should mean the existence of a morphism from $A$ to $B$.
If we specialize to the case where $A$ is the multiplicative unit 1,
so that the sequent $A\vdash B$ becomes deductively equivalent (in
linear logic) with $\vdash B$, and if we note that the unit for our
provisional $\ti$ is $(1,1,\text{true})$, then we see that truth of
$B$ should be equivalent to existence of a morphism from
$(1,1,\text{true})$ to $B$.  It is easily checked that existence of
such a morphism is precisely the second definition of truth above.

Finally, as we shall see in a moment, each definition has its own
advantages and disadvantages when one tries to prove the soundness of
linear logic, and eventually we shall need to adopt a compromise
between them.  The remark above about the relationship between
$\models_1$, $\models_2$ and negation suggests that either version of
$\models$, used alone, might have difficulties with the axioms $\vdash
A,A^\perp$ (which say that linear negation is no stronger than it
should be) or the cut rule (which says that linear negation is no
weaker than it should be).  Let us consider what happens if one tries
to establish the soundness of the axioms and cut for either version of
$\models$.  

For the axioms, we wish to show that $\ba\dpa\ba^\perp$ is true for
each object $\ba$ of \pv.  In $\ba\dpa\ba^\perp$, the questions are
pairs $(x,y)$ where $x\in A_-$ and $y\in (A^\perp)_-=A_+$, and the
answers are pairs $(a,b)$ where $a\in A_+$ and $b\in (A^\perp)_+=A_-$.
The answer $(a,b)$ is correct for the question $(x,y)$ if and only if
either $A(x,a)$ or $\neg A(b,y)$ (the latter being the definition of
$A^\perp(y,b)$).  Obviously, any question $(x,y)$ is correctly
answered by $(y,x)$.  So $\models_1\ba\dpa\ba^\perp$.  On the other
hand, we do not in general have $\models_2\ba\dpa\ba^\perp$, since an
easy calculation shows that this would mean that in $\ba$ either some
answer is correct for all questions or some question has no correct
answer.  There are, of course, easy examples of $\ba$ where this
fails; the simplest is to take $A_-=A_+=\emptyset$, and if one insists
on non-empty sets then the simplest is $A_-=A_+=\{1,2\}$ with $A$
being the relation of equality.  So, for the soundness of the axioms,
$\models_1$ works properly, but $\models_2$ does not.

Now consider the cut rule.  We wish to show that, if $\bbb\dpa\ba$ and
${\bold C}\dpa\ba^\perp$ are true, then so is $\bbb\dpa{\bold C}$.  If
we interpret truth as $\models_2$, then this is easy.  Suppose $(b,x)$
correctly answers all questions in $\bbb\dpa\ba$ and $(c,y)$ correctly
answers all questions in ${\bold C}\dpa\ba^\perp$; we claim that
$(b,c)$ correctly answers all questions $(p,q)$ in $\bbb\dpa{\bold
C}$.  Indeed, if $(p,q)$ were a counterexample, then $b$ is not
correct for $p$ and $c$ is not correct for $q$, yet $(b,x)$ is correct
for $(p,y)$ and $(c,y)$ is correct for $(q,x)$ (where the four
occurrences of ``correct'' refer to $\bbb$, $\bold C$, $\bbb\dpa\ba$,
and ${\bold C}\dpa\ba^\perp$, respectively).  But then we must have,
by definition of $\dpa$, that $x$ correctly answers $y$ in $\ba$ and
that $y$ correctly answers $x$ in $\ba^\perp$.  That is impossible, by
definition of ${}^\perp$, so the cut rule preserves $\models_2$.
Unfortunately, it fails to preserve $\models_1$.  The easiest
counterexamples occur when both $\bbb$ and $\bold C$ have questions
with no correct answers (but $B_+$ and $C_+$ are non-empty).  Then
$\bbb\dpa{\bold C}$ is not true, so the soundness of the cut rule
would require that at least one of $\bbb\dpa\ba$ and ${\bold
C}\dpa\ba^\perp$ also fail to be true.  That means that either $\ba$
or its negation must have a question with no correct answer, i.e., in
$\ba$ either some answer is correct for all questions or some question
has no correct answer.  Since that is not the case in general, we
conclude that the cut rule is unsound for $\models_1$.

Summarizing the preceding discussion, we have 
\roster
\item If we define truth allowing answers to depend on questions
($\models_1$), then the axioms of linear logic are sound but the cut
rule is not.
\item If we define truth requiring the answer to be independent of the
question ($\models_2$), then the cut rule is sound but the axioms are not.
\endroster
Fortunately, there is a way out of this dilemma.  Consider the
dependence of answers on questions that was needed to obtain the
soundness of the axioms.  At first sight, it is an extremely strong
dependence; indeed, the answer $(y,x)$ is, except for the order of
components, identical to the question $(x,y)$.  But the dependence is
special in that each component of the answer depends only on the {\sl
other\/} component of the question.  

Rather surprisingly, this sort of
cross-dependence also makes the cut rule sound.  To see this, suppose
that both $\bbb\dpa\ba$ and
${\bold C}\dpa\ba^\perp$ are true in this new sense.  That is, there
are functions $f:B_-\to A_+$ and $g:A_-\to B_+$ such that, for all
$b\in B_-$ and $x\in A_-$,
$$
B(b,g(x))\quad\text{or}\quad A(x,f(b)),
\tag1$$
and similarly there are $f':C_-\to (A^\perp)_+=A_-$ and
$g':(A^\perp)_-=A_+\to C_+$ such that, for all $c\in C_-$ and all
$y\in A_+$,
$$
C(c,g'(y))\quad\text{or}\quad\neg A(f'(c),y).
\tag2$$ 
Then we claim that $g'\circ f:B_-\to C_+$ and $g\circ f':C_-\to B_+$
satisfy, for all $b\in B_-$ and $c\in C_-$,
$$
B(b,g(f'(c))) \quad\text{or}\quad C(c,g'(f(b))),
$$
which means that $\bbb\dpa{\bold C}$ is true in the
``cross-dependence'' sense.  To verify the claim, let such $b$ and $c$
be given.  If $A(f'(c),f(b))$, then \thetag2 implies $C(c,g'(f(b)))$.
If $\neg A(f'(c),f(b))$, then \thetag1 implies $B(b,g(f'(c)))$. So the
claim is true in either case, and we have verified the soundness of
the cut rule.

By allowing the answer in one component of a sequent to depend on the
questions in the other components but not in the same component, this
``cross-dependence'' notion of truth makes crucial use of the commas
in a sequent, to distinguish the components.  But linear logic
requires (by the introduction rules for times and especially for par)
that the commas in a sequent behave exactly like the connective $\pa$.
So it seems necessary to build cross-dependence into the
interpretation of this connective.  This will lead to the correct
definition of the multiplicative connectives, replacing the
provisional interpretations given earlier.

We define the par operation on objects of \pv\ by
$$
\obj A\pa\obj B=
(A_-\times B_-, A_+^{B_-}\times B_+^{A_-}, P)
$$
where 
$$
P((x,y),(f,g))\iff A(x,f(y))\text{ or }B(y,g(x)).
$$
This operation $\pa$ is the object part of a functor, the action on
morphisms being $(f\pa g)_-=f_-\times g_-$ and $(f\pa
g)_+=f_+^{g_-}\times g_+^{f_-}$.  It is easy to check that $\pa$ is
associative (up to natural isomorphism).  In the par of several
objects, questions are tuples consisting of one question from each of
the objects, and answers are tuples of functions, each producing an
answer in one component when given as inputs questions in all the
other components.

We also interpret commas in sequents as the new $\pa$ (rather than
$\dpa$).  This change in the interpretation of the commas makes
$\models_2$ behave like the cross-dependence notion of truth described
earlier.  To see this, note that $\models_2$
requires the existence of a single answer correct for all questions at
once, but the new $\pa$ allows that answer to consist of functions
whereby each component of the answer can depend on the other
components of the question.  We therefore adopt $\models_2$ as the
(non-provisional) definition of truth, and from now on we write it
simply as $\models$.  The previous discussion shows that the axioms
and the cut rule are sound.  We sometimes refer to an answer that is
correct for all questions in an object $\ba$ as a {\sl solution\/} of
the problem $\ba$.  So truth means having a solution.

Of course, the new interpretation of par gives, by duality, a new
interpretation of times.
$$
\obj A\ti\obj B=(A_-^{B_+}\times B_-^{A_+},A_+\times B_+, T),
$$
where
$$
T((f,g),(x,y))\iff A(f(y),x)\text{ and }B(g(x),y).
$$
(This connective was called $\oslash$ in \cite{15}.)
The units for the multiplicative connectives are $1=(1,1,\text{true})$
and $\bot=(1,1,\text{false})$, where true and false represent the
obvious relations on a singleton.  The linear implication $A\imp B$
defined as $A^\perp\pa B$ is 
$$
\obj A\imp\obj B=(A_+\times B_-,A_-^{B_-}\times B_+^{A_+}, C),
$$
where 
$$
C((x,y),(f,g))\iff[A(f(y),x)\implies B(y,g(x))].
$$

Notice that a solution of $\ba\imp\bbb$ is precisely a morphism
$\ba\to\bbb$ in \pv.  This indicates that the definitions of the
multiplicative connectives and of truth, though not immediately
intuitive, are proper in the context of the category \pv.

The belief that these definitions are reasonable is reinforced by
de~Paiva's theorem \cite{15} that the multiplicative and additive
fragment of linear logic is sound for this interpretation.  (Her
theorem actually covers full linear logic, including the exponentials,
but we have not yet discussed the interpretation of the exponentials.)

Linear logic is not complete for this interpretation.  For one thing,
the interpretation validates the mix rule: If $\ba$ and $\bbb$ are
both true, then so is $\ba\pa\bbb$.  Also, the interpretation
satisfies all formulas of the form $A^\perp\pa(A\pa A)$, a special
case of weakening.  (General weakening, $A^\perp\pa(A\pa B)$ is not
satisfied; for a counterexample, take $B_+$ to be empty while all of
$A_+, A_-, B_-$ are non-empty.)

The interpretations of $\ti$ and $\pa$ remain, in spite of their
success at modeling linear logic, rather unintuitive.  This is
attested by the fact that de~Paiva \cite{15}, while using $\ti$ to
interpret the multiplicative conjunction, calls it $\oslash$ and
reserves the symbol $\ti$ for the more intuitive construction that I
called $\dti$.  \voj\ \cite{21} also discusses $\dti$, calling it
$\times$, but never has any use for $\ti$.  Since $\dti$ seems much
more natural than the ``correct'' $\ti$, it should have its own place
in the logic.

\proclaim{Suggestion 2}
Find a natural place in the theory of \pv\ for the operation $\ti$. 
\endproclaim

\head
Cardinal Characteristics of the Continuum
\endhead

We begin this section by introducing a few (just enough to serve as
examples later) of the many cardinal characteristics of the continuum
that have been studied by set-theorists, topologists, and others.  For
more information about these and other characteristics, see \cite{18} and
the references cited there.  All the cardinal characteristics
considered here (and almost all the others) are uncountable cardinals
smaller than or equal to the cardinality $\cc=2^{\aleph_0}$ of the
continuum.  So they are of little interest if the continuum hypothesis
($\cc=\aleph_1$) holds, but in the absence of the continuum hypothsis
there are many interesting connections, usually in the form of
inequalities, between various characteristics.  (There are also
independence results showing that certain inequalities are not
provable from the usual ZFC axioms of set theory.)  Part of the work
of \voj\ \cite{21, 22} on which this section is based can be viewed as a way
to extract from the inequality proofs information which is of interest
even if the continuum hypothesis holds.

\definition{Definitions}
If $X$ and $Y$ are subsets of $\nn$, we say that $X$ {\sl splits\/}
$Y$ if both $Y\cap X$ and $Y-X$ are infinite.  The {\sl splitting
number\/} $\ss$ is the smallest cardinality of any family $\Cal S$ of
subsets of $\nn$ such that every infinite subset of $\nn$ is split by
some element of $\Cal S$.  The {\sl refining number\/} (also called
the {\sl unsplitting\/} or {\sl reaping number\/}) $\rr$ is the
smallest cardinality of any family $\Cal R$ of infinite subsets of
$\nn$ such that no single set splits all the sets in $\Cal R$.  $\rs$
is the smallest cardinality of any family $\Cal R$ of infinite subsets
of $\nn$ such that, for any countably many subsets $S_k$ of $\nn$,
some set in $\Cal R$ is not split by any $S_k$.
\enddefinition

These cardinals arise naturally in analysis, for example in connection
with the Bolzano-Weierstrass theorem, which asserts that a bounded
sequence of real numbers has a convergent subsequence.  A
straightforward diagonal argument extends this to show that, for any
countably many bounded sequences of real numbers ${\bold x}_k=
(x_{kn})_{n\in\nn}$, there is a single infinite $A\subseteq\nn$ such
that the subsequences indexed by $A$, $(x_{kn})_{n\in A}$, all
converge.  If one tries to extend this to uncountably many sequences,
then the first cardinal for which the analogous result fails is $\ss$.
Also, $\rs$ is the smallest cardinality of any family $\Cal R$ of
infinite subsets of $\nn$ such that, for every bounded sequence
$(x_n)_{n\in\nn}$, there is a convergent subsequence $(x_n)_{n\in A}$
with $A\in\Cal R$.  There is an analogous description of $\rr$, where
the sequences $(x_n)_{n\in\nn}$ are required to have only finitely
many distinct terms.  For more information about these aspects of the
cardinal characteristics, see \cite{20}.

\definition{Definitions}
A function $f:\nn\to\nn$ {\sl dominates\/} another such function $g$
if, for all but finitely many $n\in\nn$, $f(n)\leq g(n)$.  The {\sl
dominating number\/} $\dd$ is the smallest cardinality of any family
$\Cal D\subseteq\fn$ such that every $g\in\fn$ is dominated by some
$f\in\Cal D$.  The {\sl bounding number\/} $\bb$ is the smallest
cardinality of any family $\Cal B\subseteq\fn$ such that no single $g$
dominates all the members of $\Cal B$.
\enddefinition

The known inequalities between these cardinals (and $\aleph_1$ and 
$\cc=2^{\aleph_0}$) are 
$$
\aleph_1\leq\ss\leq\dd\leq\cc,
$$
$$
\aleph_1\leq\bb\leq\rr\leq\rs\leq\cc,
$$
and
$$
\bb\leq\dd.
$$
It is known that any further inequalities between these cardinals are
independent of ZFC, except that it is still an open problem whether
$\rr=\rs$ is provable.

The connection between the theory of these cardinals and the category
\pv\ discussed in previous sections becomes visible when one considers
the proofs of some of these inequalities, so we shall prove the two
non-trivial (but well known) ones, $\ss\leq\dd$ and $\bb\leq\rr$.  (In
each case, only the first of the two paragraphs in the proof is
relevant to \pv, so the reader willing to take the first paragraph on
faith can skip the justification in the second paragraph.)

\demo{Proof of $\ss\leq\dd$}
There is a map $\alpha:\fn\to\pn$ sending every dominating family
$\Cal D$ (as in the definition of $\dd$) to a splitting family (as in
the definition of $\ss$).  In fact, one can associate to each infinite
$X\subseteq\nn$ a function $\beta(X)=f\in\fn$ such that, if $g$
dominates $f$, then $\alpha(g)$ splits $X$.

Given $g$, to define $\alpha(g)$, partition $\nn$ into a sequence of
intervals $[0,a_1),[a_1,a_2),\dots$ such that, for each $n\in\nn$,
$g(n)$ is at most one interval beyond $n$ (it's trivial to define such
$a_i$'s by induction), and let $\alpha(g)$ be the union of the
even-numbered intervals.  Define $\beta(X)$ to send each $n\in\nn$ to
the next element of $X$ greater than $n$.  If $f=\beta(X)$, if $g$
dominates $f$, if $a_i$'s are as in the definition of $\alpha(g)$, and
if $k$ is large enough, then the element $f(a_k-1)$ of $X$ lies in the
interval $[a_k,a_{k+1})$.  So $X$ meets all but finitely many of the
intervals $[a_k,a_{k+1})$ and is therefore split by $\alpha(g)$.
\qed\enddemo

\demo{Proof of $\bb\leq\rr$}
There is a function $\beta:\pin\to\fn$ sending every unsplittable
family $\Cal R$ (as in the definition of $\rr$) to an undominated
family (as in the definiiton of $\bb$).  In fact, one can associate to
each $g\in\fn$ a set $\alpha(g)=Y\in\pn$ such that, if $Y$ does not
split $X$ then $g$ does not dominate $\beta(X)$.

The same $\alpha$ and $\beta$ as in the preceding proof will work, as
the properties required of them here are logically equivalent to the
properties required there.
\qed\enddemo

In the notation of the preceding sections, the pair $(\beta,\alpha)$
in the first of these proofs is a morphism in \pv\ from $(\fn,\fn,
\text{is majorized by})$ to $(\pin, \pn,\text{is split by})$.  In the
second proof, we used the image of this under ${}^\perp$, namely that
$(\alpha,\beta)$ is a morphism from $(\pn,\pin,\text{does not split})$
to $(\fn,\fn,\text{does not majorize})$.  In both cases, the cardinal
inequality follows from the following general fact.  Define for each
object $\ba$ of \pv\ the {\sl norm\/} $\Vert\ba\Vert$ as the smallest
cardinality of any set $X\subseteq A_+$ of answers sufficient to
contain at least one correct answer for every question in $A_-$
(undefined if there is no such set, i.e., if some question has no
correct answer, i.e., if $\ba^\perp$ is true).  Then the existence of
a morphism $f:\ba\to\bbb$ implies that
$\Vert\ba\Vert\geq\Vert\bbb\Vert$, because $f_+$ sends any set of the
sort required in the definition of $\Vert\ba\Vert$ to one as required
for $\Vert\bbb\Vert$. (What I called the norm of $\ba$ is, in \voj's
notation \cite{21, 22} ${\frak d}(A)$; \voj's ${\frak b}(A)$ is
$\Vert\ba^\perp\Vert$.) 

It is an empirical fact that proofs of inequalities between cardinal
characteristics of the continuum usually proceed by representing the
characteristics as norms of objects in \pv\ and then exhibiting
explicit morphisms between those objects.  This fact is explicit in
\voj's \cite{21, 22} and implicit in \cite7.  It applies even to trivial
inequalities like $\bb\leq\dd$ (where the required morphism from
$(\fn, \fn, \text{is dominated by})$ to $(\fn, \fn, \text{does not
dominate})$ consists of identity maps on both components) as well as
to inequalities much deeper than the examples proved above; see for
example the presentation in \cite7 of Bartoszy\'nski's theorem \cite3
that the smallest number of meager sets whose union is not meager is
at least as large as the corresponding number for ``measure zero'' in
place of ``meager.''  

It is tempting to regard the existence of a morphism $\ba\to\bbb$ as a
strong formulation of the inequality $\Vert\ba\Vert\geq\Vert\bbb\Vert$
that is significant even in the presence of the continuum hypothesis
(which makes inequalities between cardinal characteristics trivial as
these cardinals lie between $\aleph_1$ and $\cc$ inclusive).  The
situation is, however, not quite so simple.  My student, Olga
Yiparaki, has shown that, in the presence of the continuum hypothesis
(or certain weaker assumptions), there are morphisms in \pv\ in both
directions between any two objects that correspond (as in \cite{21, 22}) to
cardinal characteristics of the continuum.  Those morphisms, however,
are highly non-constructive, whereas those used in the usual proofs of
cardinal inequalities are quite explicit.  It therefore seems likely
that a strengthening of these cardinal inequalities that retains its
significance in the presence of the continuum hypothesis is to require
not merely the existence of morphisms but the existence of ``nice''
morphisms, say ones whose components are Borel mappings.

The linear negation defined on \pv\ gives a precise version of an
intuitive ``duality'' in the theory of cardinal characteristics.  In
that theory, one often refers to the cardinals $\Vert\ba\Vert$ and
$\Vert\ba^\perp\Vert$ as being dual to each other; see for example the
introduction to \cite{13}.  On cardinals, this is not well defined, for
two objects can have the same norm while their negations have
different norms, but it is the shadow, in the world of cardinals, of
the (well defined) linear negation in \pv.  It may be worth noting in
this connection that $(\pn,\pin,\text{does not split})$, whose norm is
$\rr$, and $(\pn^\nn,\pin,\text{has no component that splits})$, whose
norm is $\rs$, have negations both of norm $\ss$.

In addition to inequalities of the sort discussed above, which relate
two cardinal characteristics of the continuum, there are a few
theorems that relate three (occasionally even four) of them.  We
consider one relatively easy example here, since it leads to an idea
that should connect to linear logic.  The example concerns Ramsey's
theorem \cite{16}, which asserts (in a simple form) that, whenever  the set
$[\nn]^2$ of two-element subsets of $\nn$ is partitioned into two
pieces, then there is an infinite $H\subseteq\nn$ that is homogeneous
in the sense that all its two element subsets lie in the same piece of
the partition.  The cardinal $\frak{hom}$ was defined in \cite5 as the
smallest cardinality of a family $\Cal H$ of infinite subsets of $\nn$
such that, for every partition of $[\nn]^2$ as in Ramsey's theorem, a
homogeneous set can be found in $\Cal H$.  It was shown in \cite5 that
this cardinal is bounded below by $\max\{\rr,\dd\}$ and above by
$\max\{\rs,\dd\}$.  The lower bound amounts to two ordinary
inequalities, $\frak{hom}\geq\rr$ and $\frak{hom}\geq\dd$, both of
which were proved by exhibiting morphisms between the appropriate
objects of \pv.  The upper bound genuinely relates three cardinals,
and we wish to make some comments about its proof, so we begin by
sketching the proof.

\demo{Proof of $\frak{hom}\leq\max\{\rs,\dd\}$}
Fix a family $\Cal R_0$ of $\rs$ subsets of $\nn$ such that no
countably many sets split all the sets in $\Cal R_0$.  Within each set
$A\in\Cal R_0$, fix a family $\Cal R_1(A)$ of $\rr$ sets such that no
single set splits them all.  Also, fix a family $\Cal D$ of functions
dominating all functions $\nn\to\nn$.  For each $A\in\Cal R_0$, for
each $B\in\Cal R_1(A)$, and for each $f\in\Cal D$, choose a subset
$Z=Z(A,B,f)$ of $B$ so thin that, if $x<y$ are in $Z$ then
$f(x)<y$.  We claim that the family $\Cal H$ of all these $Z$'s, which
clearly has cardinality $\max\{\rs,\dd\}$ (since $\rr\leq\rs$),
contains almost homogeneous sets for all partitions of $[\nn]^2$ into two
parts. ``Almost homogeneous'' means that the set becomes homogeneous
when finitely many of its elements are removed.  Since we can close
$\Cal H$ under such finite changes without increasing its cardinality,
the claim completes the proof.  

To prove the claim, let $[\nn]^2$ be partitioned into two parts.  For
each natural number $n$ let $C_n$ consist of those $x$ for which
$\{n,x\}$ is in the first part.  By choice of $\Cal R_0$, it contains
a set $A$ unsplit by any $C_n$.  Let $g(n)$ be so large that all
$x\in A$ with $x\geq g(n)$ have $\{n,x\}$ in the same piece of the
partition, and let $Q$ be the set of $n$ for which this is the first
piece.  Choose $B\in\Cal R_1(A)$ unsplit by $Q$ and $f\in \Cal D$
dominating $g$.  It is then easy to check that $Z(A,B,f)$ is almost
homogeneous for the given partition.
\qed\enddemo

To discuss this proof in terms of \pv, we introduce the natural
objects of \pv\ whose norms are the cardinals under consideration.
For mnemonic purposes, we name each object with the capital letter
corresponding to the lower-case letter naming the cardinal.
$$
\bold{HOM}=(\{p\mid p:[\nn]^2\to2\},\pin, AH),
$$
where $AH$ is the relation of almost homogeneity, $AH(p,H)$ means that
$H$ is almost homogeneous for the partition $p$.
$$
\bold D=(\fn,\fn,\text{is dominated by}).
$$
$$
\bold R=(\pn,\pin,\text{does not split}).
$$
$$
\bold R_\sigma=(\pn^\nn,\pin,\text{has no component that splits}).
$$

The structure of the preceding proof is then as follows.  From a
question $p$ in $\bold {HOM}$, we first produced a question
$(C_n)_{n\in\nn}$ in $\bold R_\sigma$.  Using an answer $A$ to this
question and also using again the original question $p$, we produced
questions $g$ in $\bold D$ and $Q$ in $\bold R$.  From answers $f$ and
$B$ to these questions, along with the previous answer $A$, we finally
produced an answer $H$ to the original question $p$ in $\bold{HOM}$.  

This can be described as a morphism into $\bold{HOM}$ from a suitable
combination of $\bold D$, $\bold R$, and $\bold R_\sigma$, but the
relevant combination is a bit different from what we have considered
previously.  The part of the construction involving $\bold D$ and
$\bold R$ is just the provisional tensor product $\dti$; that is, we
had a question $(g,Q)$ in $\bold D\dti\bold R$ and we obtained an
answer $(f,B)$ for it.  (Strictly speaking, we used a version of
$\bold R$ on $A$ rather than on $\nn$, but we shall ignore this
detail.)  The novelty is in how $\bold D\dti\bold R$ is combined with
$\bold R_\sigma$.  For what we produced from $p$ was a question in
$\bold R_\sigma$ together with a function converting answers to this
question into questions in $\bold D\dti\bold R$.  This thing that we
produced ought to be a question in the object that is being mapped to
$\bold {HOM}$.  An answer in that object ought to be what we used in
order to get the answer $H$ for $\bold{HOM}$, namely $(A,B,f)$.

Motivated by these considerations, we define a connective, denoted by
a semi-colon (to suggest sequential composition), as follows.
$$
\obj A;\obj B=(A_-\times B_-^{A_+},A_+\times B_+,S), 
$$
where 
$$
S((x,f),(a,b))\iff A(x,a)\text{ and }B(f(a),b).
$$
Thus, a question of sort $\ba;\bbb$ consists of a first question in
$\ba$, followed by a second question in $\bbb$ that may depend on the
answer to the first.  A correct answer consists of correct answers to
both of the constituent questions.  Thus, squential composition can be
viewed as describing a dialog in which the questioner first asks a
question in $\ba$, is given an answer, selects on the basis of this
answer a question in $\bbb$, and is given an answer to this as well.

The proof of $\frak{hom}\leq\max\{\rs,\dd\}$ exhibits a morphism from
$\bold R_\sigma;(\bold D\dti\bold R)$ to $\bold{HOM}$.  The cardinal
inequality follows from the existence of such a morphism, since one
easily checks that the operations on infinite cardinal norms
corresponding to the operations $\dti$ and\ ; are both simply $\max$.

The sequential composition of objects of \pv\ occurs repeatedly in the
proofs of inequalities relating three cardinal characteristics.  A
typical example is the proof that the minimum number of meager sets of
reals with a non-meager union is the minimum of $\bb$ and the minimum
number of meager sets that cover the real line \cite{13, 4, 7}.  \voj\
\cite{21} describes the strategy for proving such three-way
relations between cardinals in terms of what he calls a max-min
diagram.  This diagram amounts exacly to a morphism from the
sequential composition of two objects to a third object.  In other
words, sequential composition is the reification
of the max-min diagram  as an object of \pv.

Sequential composition also seems a natural concept to add to linear
logic from the computational point of view.  Linear logic is generally
viewed as a logic of parallel computation, but even parallel
computations often have sequential parts, so it seems reasonable to
include in the logic a way to describe sequentiality.  These ideas are
not yet sufficiently developed to support any claims about sequential
composition, as defined in the \pv\ model, being the (or a) right way
to do this.  In addition to semantical interpretations, one would
certainly want good axioms governing any sequential composition
connective that is to be added to linear logic, and one would hope
that some of the pleasant proof theory of linear logic would survive
the addition.  Much remains to be done in this direction.

\proclaim{Suggestion 3}
Find a place for sequential composition (specifically for the
connective called\ ; above) in linear logic and the theory of \pv.
\endproclaim

\head
Generalized Multiplicative Connectives
\endhead

The previous sections have led to three suggestions of natural
connectives to add to linear logic.  (Actually, Suggestion 1 concerned
not a connective but a modified notion of morphism.  But such a
modification should correspond to a reinterpretation of $\imp$ and
therefore of $\pa$ and $\ti$ as well.)  The suggested new connectives
are all analogous to the multiplicatives in that both the set of
questions and the set of answers are cartesian products.  (For the
additive connectives, one of the two sets was a disjoint union.)  The
factors in these products are either sets of questions or answers from
the constituent objects or else sets of functions, from questions to
answers or vice versa.  With these preliminary comments, it seems
natural to describe general multiplicative conjunctions (cousins of
$\ti$) as follows.

A {\sl general multiplicative conjunction\/} operates on $n$ objects
$\ba_1\dots\ba_n$ of \pv\ to produce an object $\bold C$, where 
$C_+=A_{1+}\times\dots A_{n+}$ and where $C_-$ consists of $n$-tuples
$(f_i)$ of functions where $f_i$ maps some product of $A_{j+}$'s into
$A_{i-}$. Which $A_{j+}$'s occur in the domain of which $f_i$'s is
given by the specification of the particular connective.  An answer
$(a_i)$ is correct for a question $(f_i)$ if each $a_i$ correctly
answers in $\ba_i$ the question obtained by evaluating $f_i$ at the
relevant $a_j$'s.  

For example, $\ti$ is a generalized multiplicative conjunction, for 
which $n=2$ and each $f_i$ has domain $A_j$ for the $j$ different from
$i$.  Similarly, we obtain $\dti$ if the domains of the $f_i$'s are
taken to be empty products (i.e., singletons); no $j$ is relevant to
any $i$.  Sequential composition is obtained by having $f_1$ depend on
no arguments while $f_2$ has an argument in $A_1$.

Dual (via ${}^\perp$) to generalized multiplicative conjunctions are
generalized multiplicative disjunctions.  Here the answers are allowed
to depend on some questions, rather than vice versa (exactly which
dependences are allowed is the specification of a particular
connective), and correctness means correctness in at least one
component, rather than in all.

To avoid possible confusion, we stress that the generalization of the
multiplicative connectives proposed here is quite different from that
proposed by Danos and Regnier \cite6.  The Danos-Regnier
multiplicatives can correspond to many different classical
connectives, whereas mine correspond only to conjunction and
disjunction.  One could, of course, consider combining the two
generalizations, but we do not attempt this here.

There are non-trivial {\sl unary\/} conjunction and disjunction
connectives.  The conjunction is given by 
$$
\kappa\obj A=(A_-^{A_+},A_+,\kappa A),
$$
where 
$$
\kappa A(f,a)\iff A(f(a),a).
$$
The dual disjunction is 
$$
\alpha\obj A=(A_-,A_+^{A_-},\alpha A)
$$
where
$$
\alpha A(a,f)\iff A(a,f(a)).
$$
These operations were called $T$ and $R$ in \cite{15}.

The modified concept of morphism from $\ba$ to $\bbb$ in Suggestion~1,
where $f_+$ maps $A_+\times B_-$, rather than just $B_-$, into $B_+$,
amounts to a morphism (in the standard \pv\ sense) from $\ba$ to
$\alpha\bbb$.  This concept of morphism thus gives rise to the Kleisli
category of \pv\ with respect to the monad $\alpha$.  (We have defined
$\alpha$ only on objects, but it is routine to define it on morphisms
and to describe its monad structure.)

De~Paiva's Dialectica category \cite{14} built over the category of
sets has as morphisms $\ba\to\bbb$ the \pv\ morphisms
$\kappa\ba\to\bbb$.  It is dual (via ${}^\perp$) to the category in
the preceding paragraph and is the co-Kleisli category of the comonad
$\kappa$ (see \cite{15, Prop. 7}). 

The connective $\alpha$ also provides a way to reinstate the notion of
truth $\models_1$ that was discarded when we replaced the provisional
$\ti$ and $\pa$ with the final versions.  Indeed, $\models_1\ba$ holds
if and only if $\models\alpha\ba$.

\head
Exponentials
\endhead

Girard has pointed out that the exponential connectives or modalities,\
! and \ ?, unlike the other connectives, are not determined by the
axioms of linear logic.  More precisely, if one added to linear logic
a second pair of modalities, say $!'$ and $?'$, subject to the same
rules of inference as the original pair, then one could not deduce
that the new modalities are equivalent to the old.  Several versions
of the exponentials could coexist in one model of linear logic.

\pv\ provides an example of this phenomenon.  De~Paiva \cite{15} gave an
interpretation of the exponentials in which\ ! is a combination of the
unary conjunction $\kappa$ defined above and a construction $S$ where
multisets $m$ of questions are regarded as questions and a correct
answer to $m$ is a single answer that is correct for all the questions
in $m$.  (Neither $\kappa$ nor $S$ alone can serve
as an interpretation of\ !.)  Another interpretation of the
exponentials in \pv, validating the exponential rules  of linear
logic, is given by
$$
!\obj A=(1,A_+,U)
$$
where 1 is a singleton, say $\{*\}$ and 
$$
U(*,a)\iff\forall x\in A_-\,A(x,a),
$$
and its dual
$$
?\obj A=(A_-,1,E)
$$
where 
$$
E(a,*)\iff\exists x\in A_+\,A(a,x).
$$
Intuitively, a question of type $!\ba$ (namely $*$) amounts to all
questions of type $\ba$; a correct answer in $!\ba$
must correctly answer all questions in $\ba$ simultaneously.

It is easy to check that Girard's rules of inference for the
exponentials are sound for this simple interpretation.

\enddocument